\documentclass[conference]{IEEEtran}

\IEEEoverridecommandlockouts

\def\BibTeX{{\rm B\kern-.05em{\sc i\kern-.025em b}\kern-.08em
    T\kern-.1667em\lower.7ex\hbox{E}\kern-.125emX}}

\usepackage{amsmath}
\usepackage{amssymb}
\usepackage{bm}
\usepackage{amsthm}
\usepackage{mathtools}
\usepackage{algorithm,algorithmic}
\usepackage{cleveref}
\usepackage{xcolor}
\usepackage[inline]{enumitem}
\usepackage{caption}
\usepackage{array}
\usepackage{threeparttable}
\usepackage{autonum}

\theoremstyle{definition}
\newtheorem{df}{Definition}[section]
\newtheorem{pr}[df]{Proposition}
\newtheorem{thm}[df]{Theorem}

\newtheorem{lem}[df]{Lemma}
\newtheorem{fac}[df]{Fact}
\newtheorem{rma}[df]{Remark}
\newtheorem{asm}[df]{Assumption}
\newtheorem{prob}[df]{Problem}

\newtheorem*{notation}{Notation}
\newtheorem*{relatedwork}{Related works}

\crefname{thm}{Theorem}{Theorems}
\crefname{pr}{Proposition}{Propositions}
\crefname{df}{Definition}{Definitions}
\crefname{lem}{Lemma}{Lemmas}
\crefname{fac}{Fact}{Facts}
\crefname{rma}{Remark}{Remarks}
\crefname{asm}{Assumption}{Assumptions}
\crefname{prob}{Problem}{Problems}
\crefname{ex}{Example}{Examples}

\crefname{equation}{}{}
\crefname{algorithm}{Algorithm}{Algorithms}
\crefname{figure}{Figure}{Figures}
\crefname{table}{Table}{Tables}
\crefname{chapter}{Chapter}{Chapters}
\crefname{section}{Section}{Sections}
\crefname{enumi}{}{} 

\DeclarePairedDelimiter{\abs}{\lvert}{\rvert} 
\DeclarePairedDelimiter{\norm}{\lVert}{\rVert} 
\DeclarePairedDelimiter{\rbra}{\lparen}{\rparen} 
\DeclarePairedDelimiter{\cbra}{\lbrace}{\rbrace} 
\DeclarePairedDelimiter{\sbra}{\lbrack}{\rbrack} 
\DeclarePairedDelimiter{\abra}{\langle}{\rangle} 
\DeclarePairedDelimiterX{\Set}[2]{\lbrace}{\rbrace}{#1\,\delimsize\vert\,#2}

\newcommand{\R}{\mathbb{R}}
\newcommand{\N}{\mathbb{N}}
\newcommand{\x}{\bm{x}}
\newcommand{\y}{\bm{y}}
\newcommand{\z}{\bm{z}}
\newcommand{\id}{\mathrm{Id}}

\newcommand{\frakS}{\mathfrak{S}}

\newcommand{\map}[3]{#1: #2 \to #3}
\newcommand{\seq}[2]{(#1_#2)_{#2=1}^{\infty}}
\newcommand{\del}{\partial}
\newcommand{\nb}{\nabla}
\newcommand{\D}{\mathrm{D}}
\newcommand{\inv}[1]{\frac{1}{#1}}

\newcommand{\minimize}[2]{\underset{#1\in#2}{\text{minimize }}}

\newcommand{\argmin}[1]{\underset{#1}{\text{argmin }}}
\newcommand{\mor}[2]{{}^{#2}#1}
\newcommand{\prox}[1]{\text{prox}_{#1}}
\newcommand{\toinf}{\to \infty}

\renewcommand{\baselinestretch}{0.933}

\newcommand{\noteA}{\footnote{
    Although $f$ and $g$ are weakly convex, 
    $f-g$ is actually DC function because it can be expressed as the difference of convex functions $\tilde{f} := f+\frac{\eta}{2}\norm{\cdot}^2$ and $\tilde{g} := g+\frac{\eta}{2}\norm{\cdot}^2$.
    For this reason, one might think that we could just assume $f,g$ in \cref{problem} to be convex functions instead of weakly convex functions in the first place.
    However, components of this naive DC decomposition $\tilde{f}$ and $\tilde{g}$ do not satisfy the assumption (b)(i) in \cref{problem}, namely Lipschitz continuity.
    Thus, we assume a weaker condition, i.e., weak convexity, than convexity on $f$ and $g$.}
  }
  
  \newcommand{\noteB}{\footnote{
    The functions for each entry in MCP and the capped $\ell_1$ norm is constant $\tau$ when an input is far from the origin (see \cref{tab:options of phi} in \cref{sec:introduction}).
    We set parameters of MCP and the capped $\ell_1$ norm so that $\tau=1000$.
  }}

  \newcommand{\noteD}{\footnote{
    \cref{problem} covers seemingly much more general case of the minimization of $h+\hat{f} \circ \mathfrak{S}_{1} - \hat{g}\circ\mathfrak{S}_{2}$,
    where \scalebox{0.95}[1]{$\hat{f}$}$:\mathbb{R}^{n_{1}}\to\mathbb{R}$ and $\map{\hat{g}}{\R^{n_2}}{\R}$ are Lipschitz continuous, weakly convex and prox-friendly,
    and $\map{\frakS_1}{\R^d}{\R^{n_1}}$ and $\map{\frakS_2}{\R^d}{\R^{n_2}}$ are continuously differentiable such that their Fr\'{e}chet derivative are Lipschitz continuous.
    This fact can be understood through a simple translation of this minimization into \cref{problem} by introducing $f:\R^{n_{1}+n_{2}}\to\R: $\scalebox{0.9}{$\;[\z_{1}^T, \z_{2}^{T}]^T$}$ \mapsto \hat{f}(\z_{1})$,\ $g:\R^{n_{1}+n_{2}}\to\R:$\scalebox{0.9}{$\;[\z_{1}^T, \z_{2}^T]^T$}$ \mapsto \hat{g}(\z_2)$,
    and $\mathfrak{S}:\mathbb{R}^{d}\to\mathbb{R}^{n_{1}+n_{2}}:\x\mapsto$\scalebox{0.9}{$[\mathfrak{S}_{1}(\x)^T \mathfrak{S}_{2}(\x)^T]^T$}.
  }}
  
  \newcommand{\noteE}{\footnote{
    Since $(2\eta)^{-1} \ge 1$ holds for $f$ and $g$ with parameters used in this experiment, we have $\seq{\mu}{k}\subset (0,(2\eta)^{-1}]$, and $\seq{\mu}{k}$ also satisfies \cref{eq:conditions of mu}. 
  }}
  
  \newcommand{\noteF}{\footnote{
    For a special case of \cref{problem} where $\mathfrak{S}=\id$, a variant of DCA introduced in \cite{gotoh2018dc} does not require any inner loop for the subproblem because its exact solution can be obtained by using the proximity operator of $f$.
    Moreover, another DCA-type algorithm proposed by \cite{zhang2024inexact} has a convergence guarantee even though an inexact solution of the subproblem is used, where such an inexact solution can be obtained in finite steps. 
  }}

\begin{document}

    \title{\huge VARIABLE SMOOTHING ALGORITHM FOR INNER-LOOP-FREE DC COMPOSITE OPTIMIZATIONS}
    
    \author{\IEEEauthorblockN{Kumataro Yazawa, Keita Kume, Isao Yamada\thanks{This work was supported partially by JSPS Grants-in-Aid (19H04134, 24K23885).}}
    \IEEEauthorblockA{\textit{Dept. of Information and Communications engineering, Institute of Science Tokyo}}
    Email: \{yazawa,kume,isao\}@sp.ict.e.titech.ac.jp
    }
    \maketitle

    \begin{abstract}
        We propose a variable smoothing algorithm for minimizing a nonsmooth  and nonconvex cost function. The cost function is the sum of a smooth function and a composition of a difference-of-convex (DC) function with a smooth mapping.
        At each step of our algorithm, we generate a smooth surrogate function by using the Moreau envelope of each weakly convex function in the DC function, and then perform the gradient descent update of the surrogate function.
        Unlike many existing algorithms for DC problems, the proposed algorithm does not require any inner loop.
        We also present a convergence analysis in terms of a DC critical point for the proposed algorithm
        as well as its application to robust phase retrieval.
    \end{abstract}

\section{Introduction}\label{sec:introduction}
In this paper, we consider the following nonsmooth and nonconvex optimization problem.
\begin{prob}[DC composite type problem]\label{problem}
  \begin{equation}
    \underset{\x \in \R^d}{\text{minimize }}F(\x):= h(\x)+ \underbrace{(f-g)}_{\textstyle\varphi} \circ\, \mathfrak{S}(\x), \label{eq:problem}
  \end{equation}
  where
  \begin{enumerate}[label=(\alph*)]
    \item $\map{h}{\R^d}{\R}$ is differentiable and its gradient $\map{\nb h}{\R^d}{\R^d}$ is Lipschitz continuous, i.e., there exists $L_{\nb h} >0$ such that
    $\norm{\nb h(\x) - \nb h(\y)}\le L_{\nb h} \norm{\x-\y}\ \ (\x,\y\in\R^d)$;
    \item $\map{f}{\R^n}{\R}$ and $\map{g}{\R^n}{\R}$ are
    \begin{enumerate}[label=(\roman*)]
      \item Lipschitz continuous (possibly not differentiable),
      \item weakly-convex, i.e., there exist $\eta_f,\eta_g>0$ such that $f+\frac{\eta_f}{2}\norm{\cdot}^2$ and $g+\frac{\eta_g}{2}\norm{\cdot}^2$ are convex\\(We define \scalebox{0.95}{$\eta:=\max\{\eta_f,\eta_g\}$} for convenience),
      \item prox-friendly, i.e, \textit{their proximity operators} (see \cref{moreau envelope}) are available as computable operators
    \end{enumerate}
    (see \cref{tab:options of phi} for various examples of $f-g$ in applications);
    \item $\map{\mathfrak{S}}{\R^d}{\R^n}$ is differentiable, and its Fr\'{e}chet derivative $\map{{\rm D} \mathfrak{S}}{\R^d}{\R^{n \times d}}$ is Lipschitz continuous\noteD;
    \item $F$ is bounded below, i.e., $\inf_{\x \in \R^d} F(\x)> -\infty$.
  \end{enumerate}
\end{prob}

The function $f-g$ in \cref{eq:problem} is called ``Difference-of-convex (DC) function''\noteA, and thus, we call \cref{problem} ``DC composite type problem''.

\cref{problem} appears mainly in sparsity-aware signal processing applications, 
such as image restoration \cite{you2019nonconvex}, trend filtering \cite{kim2009ell_1}, compressed sensing \cite{huang2015two}, and sparse logistic regression \cite{lu2018sparse}.
In addition, \cref{problem} also arises in robust estimation including robust phase retrieval \cite{duchi2019solving,zheng2024new}.
The DC function $f-g$ in \cref{eq:problem} is employed, e.g, to induce sparsity of the target signal translated by $\mathfrak{S}$ in sparse signal processing, or to enhance the robustness of the data fidelity to measurement outliers in robust estimation.
Such DC functions $f-g$ include 
(i) convex functions (e.g., $\ell_1$ norm \cite[Exm.24.22]{bauschke2011convex}), (ii) weakly convex functions (e.g., the minimax concave penalty (MCP) \cite{zhang2010nearly} and the smoothly clipped absolute deviation (SCAD) \cite{fan2001variable}),
and (iii) DC functions that are not weakly convex (e.g., capped $\ell_1$ norm \cite{zhang2008multi}, the trimmed $\ell_1$ norm \cite{luo2015new} and its variant \cite{sasaki2024sparse}), some of which are summarized in \cref{tab:options of phi}.
\begin{table}[t]
  \begin{threeparttable}
    
    \scriptsize
    \setlength{\tabcolsep}{0.2em} 
    \setlength{\extrarowheight}{1ex}
    \caption{
      $\varphi$ in application areas and their DC decomposition
    }
    \label{tab:options of phi}
    \begin{tabular}{|c|c|c|c|}
      \hline
      name & $\varphi(\z)=(f-g)(\z)$ & $f(\z)$ & $g(\z)$ \\ \hline
      $\ell_1$ norm & $\displaystyle\sum_{i=1}^{n}|[\z]_i|$ & $\displaystyle\sum_{i=1}^{n}|[\z]_i|$ & 0 \\ 
      MCP \cite{zhang2010nearly}& $\displaystyle\sum_{i=1}^{n} r([\z]_i)$ \tnote{*1}& $\displaystyle\sum_{i=1}^{n} r([\z]_i)$ & 0 \\ 
      Capped $\ell_1$ \cite{zhang2008multi}&
      \begin{tabular}{l}
        $\displaystyle\sum_{i=1}^{n} \min\cbra*{|[\z]_i|,\beta}$ \\[-1ex]
        (with $\beta\in\R_{++}$)
      \end{tabular}
      & $\displaystyle\sum_{i=1}^{n}|[\z]_i|$ & $\displaystyle\sum_{i=1}^{n}\max\cbra*{|[\z]_i|\text{\,--\,}\beta,0}$ \\ 
      Trimmed $\ell_1$ \cite{luo2015new}&
      \begin{tabular}{l}
        $\displaystyle\sum_{i=K+1}^{n} |[\z]_{\downarrow i}|$ \tnote{*2}\\[-1ex]
        (with $0\le K\le n-1$)
      \end{tabular}
      & $\displaystyle\sum_{i=1}^{n}|[\z]_i|$ & $\displaystyle\sum_{i=1}^{K} |[\z]_{\downarrow i}|$ \\ \hline
    \end{tabular}
  
    \vspace{0.5ex}
    \begin{tablenotes}
      \item[*1] 
      $r(t) :=
      \begin{cases*}
        \lambda|t| - \frac{t^2}{2\beta} & $|t|\le \beta \lambda$, \\
        \frac{\beta \lambda^2}{2}         & otherwise
      \end{cases*}
      $
      with $\lambda,\beta\in\R_{++}$.
      \item[*2]
      $[\z]_{\downarrow i}$ denotes the entry of $\z$ whose absolute value is the $i$-th largest.
    \end{tablenotes}
  \end{threeparttable}
  \vspace{-10pt}
\end{table}

In particular, DC functions lacking weak convexity (referred to as ``inherently DC'' functions in this paper), such as the capped $\ell_1$ norm and the trimmed $\ell_1$ norm, have been attracting great attention.
Indeed, a model using the capped $\ell_1$ norm has been reported to effectively reduce the impact of outliers in application of \textit{twin support vector machine} \cite{wang2019robust}.
Moreover, in \textit{robust principal component analysis}, which aims to decompose a matrix into a low rank matrix and a sparse matrix,
a model using the capped $\ell_1$ norm in \cite{sun2013robust} outperforms a model using $\ell_1$ norm because the capped $\ell_1$ norm avoids over-penalizing a matrix that contains entries of large absolute values, unlike the $\ell_1$ norm.
On the other hand, the trimmed $\ell_1$ norm has been utilized in terms of an exact penalty theory for cardinality-constrained optimization problems \cite{yagishita2024exactpenalizationdstationarypoints}, i.e., 
\begin{equation}\label{cardinarity constrained problem}
  \underset{\x \in \R^d }{\text{minimize }}h(\x)\quad \mathrm{s.t.} \quad \norm{\mathfrak{S}(\x)}_0 \leq K,
\end{equation}
where $\norm{\cdot}_0$ counts the number of non-zero entries of a given vector.
More precisely, if the trimmed $\ell_1$ norm is used as $f-g$,
a global (resp. local) minimizer of \cref{problem} also serves as a global (resp. local) minimizer of the cardinality-constrained problem \cref{cardinarity constrained problem} under certain conditions \cite{yagishita2024exactpenalizationdstationarypoints}.
Such an exact penalty formulation via \cref{problem} seems to be more tractable than the cardinality-constrained problem \cref{cardinarity constrained problem} from a viewpoint of designing algorithms.

A commonly used existing approach to minimization of DC functions is DC algorithm (DCA) (see, e.g, \cite{le2024open}).
At each iteration of DCA, a subtrahend convex function in a DC function is replaced with an affine minorization by utilizing its subgradient, and then minimizes the resulting surrogate function as a subproblem.
For \cref{problem}, \textit{DC composite algorithm} (DCCA) \cite{le2024minimizing}, which is an extension of DCA, can be employed.
If an exact solution to the subproblem in DCCA is available, then DCCA has a convergence guarantee in terms of \textit{DC critical point} (see \cref{df:critical point}).
In practice, however, DCCA requires infinite iterations of an inner loop so as to find the exact solution of the subproblem\noteF. 
The convergence analysis of DCCA does not cover realistic cases where only inexact solutions of the subproblem are available.

In this paper, we propose a variable smoothing algorithm for \cref{problem} that does not require any inner loop for the subproblem.
Our algorithm (\cref{algorithm}) is designed as a gradient descent update of a time-varying smoothed surrogate function of $F$ in \cref{eq:problem}.
With \textit{the Moreau envelopes} (see \cref{moreau envelope}) $\mor{f}{\mu}$ of $f$ and $\mor{g}{\mu}$ of $g$, the proposed surrogate function is given as $h + (\mor{f}{\mu_k}-\mor{g}{\mu_k})\circ \mathfrak{S}$, 
where $\seq{\mu}{k} \subset \R$ is a monotonically decreasing sequence of convergence to zero.
We also present an asymptotic convergence analysis of the proposed algorithm in the sense of a DC critical point (see \cref{thm:convergence theorem,rm:convergence analysis}).
To verify effectiveness of  the proposed model (i.e, \cref{problem}) and the proposed algorithm, we conduct numerical experiments in a scenario of the robust phase retrieval (e.g., \cite{duchi2019solving,zheng2024new}) with its new optimization model.
\begin{relatedwork}
Our algorithm serves as an extension of algorithms \cite{bohm2021variable,kume2024variable,kume2024variableLong} proposed for \cref{problem} in a special case where $g\equiv0$ (more precisely, $\mathfrak{S}$ is linear in \cite{bohm2021variable}).
The existing algorithms in \cite{bohm2021variable,kume2024variable,kume2024variableLong} can be applied only to optimization problems involving weakly convex functions, 
while the proposed algorithm can cover even inherently DC functions. 
\end{relatedwork}

\begin{notation}
  $\N$, $\R$ and $\R_{++}$ denote respectively the sets of all positive integers, all real numbers and all positive real number.
  $\norm{\cdot}$ and $\abra{\cdot,\cdot}$ are respectively the Euclidean norm and the standard inner product.
  For $\bm{v}\in\R^n$, $[\bm{v}]_i\in\R$ stands for the $i$-th entry.
  We use $\mathrm{Id}$ to denote the identity mapping.
  For Euclidean spaces $\mathcal{X},\mathcal{Y}$ and a continuously differentiable mapping $J:\mathcal{X} \to \mathcal{Y}$, its Fr\'{e}chet derivative at $\x \in \mathcal{X}$ is 
  the linear operator $\map{\mathrm{D}J(\x)}{\mathcal{X}}{\mathcal{Y}}$ such that $\lim_{\mathcal{X}\setminus{\{\bm{0}\}}\ni\bm{h} \to \bm{0}}\frac{J(\x+\bm{h})-J(\x)-DJ(\x)[\bm{h}]}{\norm{\bm{h}}} = 0$.
  In particular with $\mathcal{Y}=\R$, $\map{\nb J}{\mathcal{X}}{\mathcal{X}}$ is called the gradient of $J$ if
  $\nb J(\x) \in \mathcal{X}$ at $\x \in \mathcal{X}$ satisfies $\mathrm{D}J(\x)[\bm{v}] = \abra{\nb J(\x),\bm{v}}\ (\bm{v}\in\mathcal{X})$.
\end{notation}

\section{Preliminary}
As an extension of the subdifferential of convex functions, we use the following subdifferential of nonconvex functions.
(see, e.g., a recent survey \cite{li2020understanding} for readers who are unfamiliar with nonsmooth analysis).
\begin{df}[Regular subdifferential {\cite[Def. 8.3]{rockafellar2009variational}}]
  For a function $\map{\phi}{\R^d}{\R}$, \textit{the regular subdifferential} of $\phi$ at $\bar{\x} \in \R^d$, denoted as $\del \phi(\bar{\x}) \subset \R^d$, is the set of all vectors $\bm{v} \in \R^n$ such that
  \begin{equation}
    \lim_{\delta \searrow 0}\ \inf_{0 < \norm{\x - \bar{\x}} < \delta}\frac{\phi(\x)-\phi(\bar{\x})-\abra{\bm{v},\x-\bar{\x}}}{\norm{\x-\bar{\x}}} \ge 0.
  \end{equation}
\end{df}

If $\phi$ is convex, this regular subdifferential is equivalent to the convex subdifferential \cite[Proposition 8.12]{rockafellar2009variational}.
Furthermore, if $\phi$ is Fr\'{e}chet differentiable at $\bar{\x}$, $\del\phi(\bar{\x}) = \cbra*{\nb \phi(\bar{\x})}$ holds \cite[Exercise 8.8(a)]{rockafellar2009variational}.

Unfortunately, finding a global minimizer of \cref{problem} is not realistic due to the severe nonconvexity of $F$.
Instead, in this paper, we focus on finding a DC critical point defined, with the regular subdifferentials, as follows.
\begin{df}[DC critical point for \cref{problem} \cite{le2024minimizing}]\label{df:critical point}
  A point $\x^\star \in \R^d$ is said to be a \textit{DC critical point} for \cref{problem} if
  \begin{equation} \label{eq:critical point}
    \del (h+f \circ \mathfrak{S})(\x^\star)- \del (g \circ \mathfrak{S})(\x^\star) \ni \bm{0}.
  \end{equation}
\end{df}

\begin{lem}[Relationship between local minimizer and DC critical point]\label{critical point is minimizer}
Let  $\x^\star \in \R^d$ be a local minimizer of $F$ in \cref{problem}. Then, $\x^\star$ is a DC critical point for \cref{problem}.
\end{lem}

From \cref{critical point is minimizer}, being a DC critical point is a necessary condition for being a local minimizer.
Moreover, finding such a DC critical point has been used as an acceptable goal in many DC optimization literature \cite{le2024open,le2024minimizing,gotoh2018dc,zhang2024inexact,sun2023algorithms}.

The Moreau envelope plays an important role in this paper for designing the proposed algorithm. 
\begin{df}[Moreau envelope, proximity operator \cite{bohm2021variable}] \label{moreau envelope}
  Let $\map{\psi}{\R^n}{\R}$ be an $\eta_\psi$-weakly convex function with $\eta_\psi > 0$. Its Moreau envelope and proximity operator at $\bar{\z} \in \R^n$ with $\mu \in (0,\eta_\psi^{-1})$ are respectively defined as
  \begin{align}
    \mor{\psi}{\mu}(\bar{\z}) &:= \min_{\z \in \R^n} \cbra*{\psi(\z) + \frac{1}{2\mu}\norm*{\z-\bar{\z}}^2}, \\
    \prox{\mu \psi}(\bar{\z})  &:= \argmin{\z \in \R^n} \cbra*{\psi(\z) + \frac{1}{2\mu}\norm*{\z-\bar{\z}}^2},
  \end{align}
  where $\prox{\mu \psi}$ is single-valued due to the strong convexity of $\psi + (2\mu)^{-1}\norm{\cdot - \bar{\z}}^2$.
\end{df}
The Moreau envelope $\mor{\psi}{\mu}$ serves as a smoothed surrogate function of $\psi$
because of the next properties.
\begin{fac}[Properties of Moreau envelope]\label{Properties of Moreau envelope}
  Let $\map{\psi}{\R^n}{\R}$ be an $\eta_\psi$-weakly convex function with $\eta_\psi > 0$.
  For $\mu \in (0,\eta_\psi^{-1})$, the following hold.
  \begin{enumerate}[label=(\alph*),leftmargin=2em]
    \item \cite[Theorem 1.25]{rockafellar2009variational} $(\z \in \R^n)\ \lim_{\mu \searrow 0} \mor{\psi}{\mu}(\bm{z}) = \psi(\bm{z})$.
    \item \cite[Collorary 3.4]{hoheisel2020regularization} $\mor{\psi}{\mu}$ is continuously differentiable with $(\z \in \R^n)\ \nb\mor{\psi}{\mu}(\z) = \mu^{-1}\rbra*{\z - \prox{\mu \psi}(\z)}$.
    \item \cite[Collorary 3.4]{hoheisel2020regularization} $\nb\mor{\psi}{\mu}$ is Lipschitz continuous with $L_{\nb\mor{\psi}{\mu}}:=\max\cbra{\mu^{-1},\frac{\eta_\psi}{1-\eta_\psi\mu} }$.
  \end{enumerate}
\end{fac}
Note that for $f$ and $g$ in \cref{problem}, we can compute $\nb\mor{f}{\mu}$ and $\nb\mor{g}{\mu}$ in closed forms because these functions are assumed to be prox-friendly (see the assumption (b)(iii) of \cref{problem}). 

\section{Variable Smoothing Algorithm for \\ DC Composite Problem}
\subsection{Design of Smooth Surrogate Function}
In our algorithm, we use the following function as a smooth surrogate function of $F$ in place of the direct utilization of the nonsmooth function $F$.
\begin{df}[Surrogate function] \label{df:surrogate function}
  Consider \cref{problem}.
  For $\mu \in \rbra{0,\eta^{-1}}$, we define a surrogate function of the cost function $F$ in \cref{problem} as
  \begin{equation}
    F^{\abra{\mu}} := h + (\mor{f}{\mu}-\mor{g}{\mu})\circ \mathfrak{S} \label{eq:surrogate function}.
  \end{equation}
\end{df}
By \cref{Properties of Moreau envelope}(b), $F^{\abra{\mu}}$ is also continuously differentiable.

The next theorem suggests how to find a DC critical point in \cref{eq:critical point} using the surrogate function $F^{\abra{\mu}}$.
\begin{thm}\label{thm:bound of distance}
  Consider \cref{problem}.
  Suppose that a positive sequence $\seq{\mu}{k} \subset \rbra{0,\eta^{-1}} $ converges to 0.
  For the function sequence $F_k := F^{\abra{\mu_k}}\ (k\in\N)$ with \cref{eq:surrogate function} and any convergent sequence $\seq{\x}{k} \subset \R^d \to \exists \bar{\x} \in \R^d$, we have
  \begin{equation}\label{eq:bound of distance}
    \text{dist}\rbra{\bm{0},\del(h+f\circ \mathfrak{S})(\bar{\x})-\del(g \circ \mathfrak{S})(\bar{\x})} \le \liminf_{k \to \infty}\norm{\nb F_k(\x_k)},
  \end{equation}
  where $\text{dist}(\bm{v},S):=\inf_{\bm{w}\in S}\norm{\bm{v}-\bm{w}}$ for a point $\bm{v}\in\R^d$ and a set $S\subset\R^d$.
\end{thm}
\cref{thm:bound of distance} implies that $\bar{\x}$ is a DC critical point in the sense of \cref{eq:critical point} if the right hand side of \cref{eq:bound of distance} is zero.
Hence, our goal of finding a DC critical point of \cref{problem} is reduced to designing an algorithm to generate a point sequence $\seq{\x}{k}$ such that $\liminf_{k \to \infty}\norm{\nb F_k(\x_k)}=0$.
\subsection{Proposed Algorithm and Its Convergence Analysis}
We propose \cref{algorithm} based on the gradient descent method of the smoothed surrogate function $F_k:=F^{\abra{\mu_k}}$.
\begin{algorithm}[t]
  \caption{Variable smoothing algorithm for DC composite type problem (\cref{problem})}
  \label{algorithm}
  \begin{algorithmic}[1]
      \REQUIRE $\x_1 \in \R^d,\seq{\mu}{k}\subset(0,(2\eta)^{-1}]$ enjoying \cref{eq:conditions of mu}.
      \FOR {$k=1,2,3,\dots$}
      \STATE Set $F_k:= F^{\abra{\mu_k}}=h + \rbra*{\mor{f}{\mu_k} - \mor{g}{\mu_k}} \circ \mathfrak{S}$
      \STATE Obtain $\gamma_k$ by \cref{backtracking}
      \STATE $\x_{k+1} \leftarrow \x_k - \gamma_k\nb F_k(\x_k)$
      \ENDFOR
  \end{algorithmic}
\end{algorithm}

We design $\seq{\mu}{k}\subset(0,(2\eta)^{-1}]$ to satisfy the following condition (introduced in \cite{kume2024variableLong}) so as to establish a convergence analysis of \cref{algorithm}:
\begin{equation}\label{eq:conditions of mu}
  \left\{
    \begin{aligned}
      &\text{(i) }\textstyle\lim_{k\toinf}\mu_k = 0, \quad \text{(ii) }\textstyle\sum\nolimits_{k=1}^{\infty}\mu_k = \infty,\\
      &\text{(iii) }(\exists M \ge 1, \forall k \in \N)\ 1 \le \mu_k/\mu_{k+1} \le M.
    \end{aligned}
  \right.
\end{equation}
For example, $\mu_k:=(2\eta)^{-1}k^{-\inv{\alpha}}$ with $\alpha\ge1$ enjoys the condition \cref{eq:conditions of mu}
($\alpha=3$ is reported to be an appropriate value for a reasonable convergence rate of a special case of \cref{algorithm} with $g\equiv 0$ \cite{bohm2021variable,kume2024variableLong}).

To obtain a stepsize $\gamma_k$ in line 3 of \cref{algorithm}, we employ the so-called \textit{backtracking algorithm} in \cref{backtracking} which has been utilized as a standard stepsize selection for smooth optimization (see, e.g., \cite{andrei2020nonlinear}).
\begin{algorithm}[t]
  \caption{Backtracking algorithm to find $\gamma_k$}
  \label{backtracking}
  \begin{algorithmic}[1]
    \REQUIRE $\gamma_\text{initial}>0,\ \rho\in(0,1),\ c\in(0,1)$
    \STATE $\gamma_k\leftarrow\gamma_\text{initial}$
    \WHILE{\scalebox{0.9}{$F_k(\x_k-\gamma_k \nb F_k(\x_k)) > F_k(\x_k) - c\gamma_k\norm{\nb F_k(\x_k)}^2$}}
    \STATE $\gamma_k \leftarrow \rho\gamma_k$
    \ENDWHILE
    \ENSURE $\gamma_k$
  \end{algorithmic}
\end{algorithm}
The finite termination of \cref{backtracking} is guaranteed under the following assumption.
\begin{asm}[Decent assumption]\label{asm:decent lemma}
  Consider the surrogate function $F^{\abra{\mu}}$ in \cref{eq:surrogate function} with $\mu\in(0,(2\eta)^{-1}]$.
  Assume that for any $\x,\y\in\R^d$,
  \begin{equation}\label{eq:decent lemma}
    F^{\abra{\mu}}(\y) \le F^{\abra{\mu}}(\x) + \abra{\nb F^{\abra{\mu}}(\x), \y-\x} + \frac{\kappa_\mu}{2}\norm{\y-\x}^2,
  \end{equation}
  where $\kappa_\mu = \varpi_1 + \frac{\varpi_2}{\mu}$ with some $\varpi_1,\varpi_2\in\R_{++}$.
\end{asm}
\begin{rma}[Sufficient conditions of \cref{asm:decent lemma}]
  In analogy with \cite[Lemma 3.2]{kume2024variable}, if $\frakS$ is Lipschitz continuous, then $\nb F^{\abra{\mu}}$ turns out to be Lipschitz continuous with a Lipschitz constant
  \begin{equation}
    L_{\nb F^{\abra{\mu}}} = L_{\nb h} + L_{\D\mathfrak{S}}\rbra*{L_f+L_g} + \frac{2L_{\mathfrak{S}}^2}{\mu}, \label{eq:simple Lipschitz constant}
  \end{equation}
  where each $L_\Theta$ denotes a Lipschitz constant of a mapping $\Theta$.
  Because $\kappa_\mu = L_{\nb F^{\abra{\mu}}}$ enjoys \cref{eq:decent lemma} \cite[Lemma 5.7]{beck2017first}, \cref{asm:decent lemma} is achieved in this case.
  On the other hand, we also will present, in Section IV, an example where Assumption III.3 is satisfied without the Lipchitz continuity of $\frakS$ (see \cref{pr:example satisfying assumption}).
\end{rma}

Here, in order to show $\liminf_{k \to \infty}\norm{\nb F_k(\x_k)}=0$ with $\seq{\x}{k}$ generated by \cref{algorithm}, we present the following lemma to see a behavior of the gradient sequence $\rbra*{\nb F_k(\x_k)}_{k=1}^{\infty}$.
\begin{lem}\label{thm:lemma for convergence analysis}
  Consider \cref{problem}. Choose arbitrarily a sequence $\seq{\mu}{k} \subset (0,(2\eta)^{-1}]$ satisfying \cref{eq:conditions of mu}, an initial point $\x_1 \in \R^{d}$, and inputs of \cref{backtracking} $(\gamma_\text{initial},\rho, c)\in\R_{++}\times(0,1)\times(0,1)$.
  Under \cref{asm:decent lemma}, the following inequality holds
  for the function sequence $F_k$ and the point sequence $\seq{\x}{k}$ produced by \cref{algorithm}:
  \begin{equation}
    (\underline{k},\bar{k} \in \N\ \text{s.t.}\ \underline{k}\le\bar{k})\ \min_{\underline{k}\le k \le \bar{k}}\norm{\nb F_k(\x_k)} \le \sqrt{\frac{C}{\sum_{k=\underline{k}}^{\bar{k}}\mu_k}},\label{eq:lemma for convergence analysis}
  \end{equation}
  where $C>0$ is a constant.
\end{lem}
By noting that $\sum_{k=1}^{\infty}\mu_k = \infty$ from the condition (ii) in \cref{eq:conditions of mu}, we obtain the next convergence theorem.
\begin{thm}\label{thm:convergence theorem}
  Under the setting of \cref{thm:lemma for convergence analysis}, we have
  \begin{equation}\label{eq:convergence theorem}
      \liminf_{k \to \infty}\norm{\nb F_k(\x_k)} = 0. 
  \end{equation}
\end{thm}

\begin{rma}[Interpretation of \cref{thm:convergence theorem}]\label{rm:convergence analysis}
  \cref{thm:convergence theorem} means that we can choose a subsequence $(\x_{m(l)})_{l=1}^{\infty}$ such that $\lim_{l \to \infty}\norm{\nb F_{m(l)}(\x_{m(l)})} = 0$,
  where $\map{m}{\N}{\N}$ is monotonically increasing.
  Every cluster point of $(\x_{m(l)})_{l=1}^{\infty}$ is guaranteed to be a DC critical point of Problem 1.1 by applying \cref{thm:bound of distance}.
\end{rma}

\section{Application to robust phase retrieval}\label{sec:experiment}
\subsection{Optimization model in robust phase retrieval}
To demonstrate the effectiveness of the proposed model (\cref{problem}) and \cref{algorithm},
we carried out numerical experiments in a scenario of \textit{the phase retrieval}.
The phase retrieval is widely used, e.g., for crystallography \cite{millane1990phase}, optical imaging \cite{shechtman2015phase}, and astronomy \cite{fienup1987phase}.
The phase retrieval is a problem of estimating an original signal $\x^\star$ or $-\x^\star \in \R^d$ from the magnitude measurement
\begin{equation}
  \scalebox{0.95}{\text{$\bm{b}^\star := (A\x^\star) \odot (A\x^\star) := \sbra*{\abra{\bm{a}_1,\x^\star}^2,\abra{\bm{a}_2,\x^\star}^2,...,\abra{\bm{a}_n,\x^\star}^2}^T$}}
\end{equation}
where $A = [\bm{a}_1^T,\bm{a}_2^T...,\bm{a}_n^T]^T\in\R^{n\times d}$, and $\odot$ means the element-wise product.
While this $\bm{b}^\star$ is clean measurement, a measurement in real-world applications may be corrupted by a noise.
In particular, \cite{duchi2019solving} considers the measurements $\bm{b}\in\R^n$ with outliers as:
\begin{equation}\label{eq:corrrupted measurement}
  [\bm{b}]_i := 
  \begin{cases}
    \abra{\bm{a}_i,\x^\star}^2 & i\in \mathcal{I}_{\text{in}}\\
    \xi_i & i\in \mathcal{I}_{\text{out}}
  \end{cases}
\end{equation}
where $\mathcal{I}_{\text{in}},\mathcal{I}_{\text{out}}\subset \{1,2,...,n\}$ denote disjoint index sets of inliers and outliers such that $\mathcal{I}_\text{in} \cup \mathcal{I}_\text{out} = \{1,2,\ldots,n\}$,
and $\xi_i > 0$ is a random noise.

In order to circumvent performance degradation caused by outliers, a robust phase retrieval \cite{duchi2019solving} has been formulated as
\begin{equation}\label{eq:model with l1}
  \minimize{\x}{\R^d} \norm*{(A\x) \odot (A\x)-\bm{b}}_1
\end{equation}
with the $\ell_1$ norm $\map{\norm{\cdot}_1}{\R^n}{\R},\ \z \mapsto  \sum_{i=1}^{n} |[\z]_i|$.
Although good estimation results have been reported by solving \cref{eq:model with l1} with \textit{Proximal linear algorithm} \cite{duchi2019solving} and \textit{Inexact proximal linear algorithm} \cite{zheng2024new},
it is questionable whether the $\ell_1$ norm in \cref{eq:model with l1} can adequately suppress the effects of the outliers. 
To explain this, we rewrite the cost function in \cref{eq:model with l1} as
$\sum_{i\in\mathcal{I}_{\text{in}}}\abs*{\abra{\bm{a}_i,\x}^2-\abra{\bm{a}_i,\x^\star}^2} + 
\sum_{i\in\mathcal{I}_{\text{out}}}\abs*{\abra{\bm{a}_i,\x}^2-\xi_i}$.
If the cardinality $\# \mathcal{I}_\text{out}$ and each $\xi_i$ are large, then the second summation also becomes large even if $\x$ is close to $\x^\star$ or $-\x^\star$.
Such a situation may lead to performance degradation.

To resolve this issue, we propose the following reformulation of the robust phase retrieval.
\begin{prob}[Proposed model for robust phase retrieval]\label{prob:model with nonconvex}
  For given $A\in\R^{n\times d}$ and $\bm{b}\in\R^n$,
  \begin{equation}\label{eq:model with nonconvex}
    \minimize{\x}{\R^d} F(\x) := \varphi \rbra*{(A\x) \odot (A\x)-\bm{b}},
  \end{equation}
  where $\varphi$ is chosen from \cref{tab:options of phi} in \cref{sec:introduction}.\\
  (Note: this problem is a special case of \cref{problem}, where $h\equiv 0$, $f-g = \varphi$, and $\frakS : \x \mapsto \rbra*{A\x}\odot\rbra*{A\x}-\bm{b}$.)
\end{prob}

The proposed model \cref{eq:model with nonconvex} with nonconvex $\varphi$ such as MCP \cite{zhang2010nearly}, the capped $\ell_1$ norm \cite{zhang2008multi}, and the trimmed $\ell_1$ norm \cite{luo2015new} seems to be more robust for the large outliers than the existing model \cref{eq:model with l1}.
Indeed, the function for each entry in MCP and the capped $\ell_1$ does not exceed a certain tunable constant value even if its entry has a large absolute value, while $\ell_1$ norm does not have such desirable property.
Alternatively, a large entry tends to be excluded from the summation in the trimmed $\ell_1$ norm because the $K$ largest absolute values are ignored therein.
Therefore, these nonconvex $\varphi$ are expected to remedy over-penalization in the model \cref{eq:model with l1}.

\cref{algorithm} can be employed for \cref{prob:model with nonconvex} since it is the special case of \cref{problem}.
Furthermore, the proposed convergence analysis in \cref{thm:convergence theorem} can be applied to \cref{prob:model with nonconvex} because \cref{asm:decent lemma} is achieved as shown in the following proposition.
To the best of the authors' knowledge, \cref{algorithm} is the first inner-loop-free algorithm applicable to \cref{prob:model with nonconvex} with nonconvex $\varphi$.
\begin{pr}\label{pr:example satisfying assumption}
  For \cref{prob:model with nonconvex} with $\varphi$ chosen from Table 1, $F^{\abra{\mu}}$ in \cref{eq:surrogate function} with $\mu\in(0,(2\eta)^{-1}]$ satisfies \cref{asm:decent lemma}.
\end{pr}

\subsection{Numerical experiments}
In order to evaluate estimation performance of the robust phase retrieval via the proposed model \cref{eq:model with nonconvex},
we applied \cref{algorithm} to the model \cref{eq:model with nonconvex} using $\ell_1$ norm, MCP, the capped $\ell_1$ norm, and the trimmed $\ell_1$ norm as $\varphi$ (see \cref{tab:options of phi} in \cref{sec:introduction}).
Note that \cref{eq:model with nonconvex} with $\ell_1$ norm is the same as the existing model \cref{eq:model with l1}, e.g., in \cite{duchi2019solving,zheng2024new}.

Our experimental setup inspired by \cite{zheng2024new} is as follows.
We drew each entry of $A\in\R^{200\times50}$ from the normal distribution $\mathcal{N}(0,1)$.
Each entry of the original signal $\x^\star\in\R^{50}$ was chosen from 1 or -1 with a probability of 0.5 respectively.
The number of outliers $\mathcal{I}_\text{out}$ was 10, that was 5\% of all entries of $\bm{b}\in\R^{200}$, and the position of outliers was randomly chosen.
The value of each outliers $\xi_i$ was given by
\begin{equation}\label{eq:choice of xi}
  \xi_i = \Omega \tan\rbra*{\frac{\pi}{2} u_i}\ (\ge 0),
\end{equation}
where $u_i$ was drawn from the uniform distribution of $[0,1]$, and $\Omega>0$ was used to control the magnitude of $\xi_i$.
We used the parameter $(\lambda,\beta)=(1,2000)$ and $(2,500)$ for MCP, $\beta = 1000$ for the capped $\ell_1$ norm\noteB, and $K = 10\,(=\#\mathcal{I}_\text{out})$ and $20$ for the trimmed $\ell_1$ norm.
In \cref{algorithm}, we employed $\x_1 \sim \mathcal{N}(\bm{0}, I_{50})$ for a random initial point, 
$\mu_k = k^{-\inv{3}}\ (k\in\N)$ for parameters of the Moreau envelope\noteE, and $(\gamma_\text{initial},\rho,c) = (1,0.8,0.0001)$ for inputs of \cref{backtracking}.
We stopped \cref{algorithm} when $\norm{\nb F_k(\x_k)} < 0.001$ held or the iteration $k$ reached to 10000.
For each $\varphi$ in \cref{tab:options of phi} and each $\Omega\in\{100,1000,2000,5000,10000\}$ in \cref{eq:choice of xi}, we performed 50 trials of estimation with the model \cref{eq:model with nonconvex}.
We judged that an estimation succeeds if the relative error, defined as $\min\cbra*{\norm{\x^\star-\x^\diamond },\norm{\x^\star+\x^\diamond}}/\norm{\x^\star}$, achieves a smaller value than $10^{-3}$, where $\x^\diamond$ is the final estimate of \cref{algorithm}.
As in \cite{zheng2024new}, we used an estimation performance criterion called ``success rate'' that is the percentage of the successful estimation out of 50 estimations.

\cref{tab:success rate} shows the success rates for each $\varphi$ and each $\Omega$.
From \cref{tab:success rate}, the model \cref{eq:model with nonconvex} with nonconvex functions, i.e., MCP, the capped $\ell_1$ norm, and the trimmed $\ell_1$ norm, keep high success rates even with large outliers while the model \cref{eq:model with nonconvex} with the $\ell_1$ norm does not.
In particular, the inherently DC functions, the capped $\ell_1$ norm and the trimmed $\ell_1$ norm with $K=10$, have higher success rates than others when outliers are large. 
(Note that the result of the trimmed $\ell_1$ with $K=10$ is based on the utilization of the number of outliers.)
\cref{tab:time} demonstrates the averaged CPU time for \cref{algorithm} to be terminated.
\cref{algorithm} for the model \cref{eq:model with nonconvex} achieves the fastest average convergence speed (i) with the capped $\ell_1$ among all estimations for huge $\Omega \in \{5000, 10000\}$, and (ii) with the trimmed $\ell_1$ among successful estimations for all $\Omega$. 
As above, the proposed model using inherently DC functions $\varphi$ have better results in both estimation performance and speed than the existing model using $\ell_1$ norm.
\begin{table}[t]
  \scriptsize
  \setlength{\tabcolsep}{0.3em}
  \setlength{\extrarowheight}{0.1ex}
  \centering
  \caption{Success rate [\%]}
  \begin{tabular}{|l||r|r|r|r|r|r|}
    \hline
     & Existing &\multicolumn{5}{c|}{Proposed} \\\hline
    $\Omega$ & $\ell_1$ & 
    \begin{tabular}{c}
      MCP \\($\lambda=1$)
    \end{tabular} & 
    \begin{tabular}{c}
      MCP \\($\lambda=2$)
    \end{tabular} & Capped $\ell_1$ & 
    \begin{tabular}{c}
      Trimmed $\ell_1$\\($K=10$)
    \end{tabular} &
    \begin{tabular}{c}
      Trimmed $\ell_1$ \\($K=20$)
    \end{tabular} 
     \\ \hline \hline
    10 &  90 & 86 & 88 & \textbf{92} & 66 & 54\\ \hline
    1000 & 60 & \textbf{80} & 78 & 70& 78 & 72\\ \hline
    3000 & 56 & \textbf{84} & 80 & \textbf{84}& 82 & 74\\ \hline
    5000 & 54 & 80 & 82 & \textbf{86} & 84 & 78\\ \hline
    10000 & 54 & 80 & 82 & 86 & \textbf{88} & 78\\ \hline
  \end{tabular}
  \label{tab:success rate}
\end{table}
\begin{table}[t]
\begin{threeparttable}
  \setlength{\tabcolsep}{0.26em} 
  \setlength{\extrarowheight}{0.1ex}
  \centering
  \caption{Averaged time [sec]}
  \fontsize{6.8pt}{9pt}\selectfont
  \begin{tabular}{|l||c|c|c|c|c|c|}
    \hline
    & Existing &\multicolumn{5}{c|}{Proposed} \\\hline
    $\Omega$ & $\ell_1$ & 
    \begin{tabular}{c}
      MCP \\($\lambda=1$)
    \end{tabular} & 
    \begin{tabular}{c}
      MCP \\($\lambda=2$)
    \end{tabular} & Capped $\ell_1$ & 
    \begin{tabular}{c}
      Trimmed $\ell_1$\\($K=10$)
    \end{tabular} &
    \begin{tabular}{c}
      Trimmed $\ell_1$ \\($K=20$)
    \end{tabular}  \\ \hline \hline
    10 & \textbf{2.50} (2.50)& 3.20 (3.22) & 3.30 (3.32) & 3.30 (3.31)& 5.74 (\textbf{0.69}) & 8.26 (1.52)\\ \hline
    1000 & \textbf{2.55} (2.57)& 3.20 (3.22)& 3.26 (3.26) & 3.32 (3.35)& 3.74 (\textbf{0.55})& 4.95 (0.71)\\ \hline
    3000 & \textbf{2.47} (2.50)& 3.23 (3.25) & 3.22 (3.23) & 2.83 (2.75)& 3.62 (0.76)& 4.48 (\textbf{0.52})\\ \hline
    5000 & 2.55 (2.61)& 3.11 (3.12) & 3.16 (3.17) & \textbf{2.36} (2.21)& 3.12 (0.76)& 3.61 (\textbf{0.46})\\ \hline
    10000 & 2.49 (2.55) & 3.13 (3.14) & 3.08 (3.08) & \textbf{1.65} (1.39) & 2.50 (0.75) & 3.74 (\textbf{0.48})\\ \hline
  \end{tabular}
  \label{tab:time}
  \begin{tablenotes}
    \item[] The values out of and in the parentheses are the averaged time taken for all estimations and only for successful estimations, respectively.
  \end{tablenotes}
\end{threeparttable}
\vspace{-10pt}
\end{table}

\section{Conclusion}
We presented an inner-loop-free variable smoothing algorithm for nonsmooth DC composite type problems with its convergence analysis.
The proposed algorithm was designed to find a DC critical point by generating the sequence of points at which the gradient of the smooth surrogate function approaches zero.
The numerical experiments in a scenario of robust phase retrieval demonstrated the effectiveness of the proposed optimization model \cref{eq:model with nonconvex} using inherently DC functions.

\bibliographystyle{IEEEtran}
\bibliography{IEEEabrv,my_bib}

\end{document}